\documentclass[reqno, 11pt]{amsart} \usepackage{amsmath,amsthm,amscd,amsfonts,amssymb}
\usepackage{latexsym}
\usepackage[usenames]{color}
\newcommand{\CC}{\mathbb C}

\newcommand{\RR}{\mathbb R}

\newcommand{\ZZ}{\mathbb Z}
\newcommand{\EE}{\mathbb E}

\newcommand{\nn}{\mathcal N}

\newcommand{\ml}{\mathcal L}

\newtheorem{thm}{Theorem}[section]
\newtheorem{lemma}[thm]{Lemma}
\newtheorem{cor}[thm]{Corollary}
\newtheorem{prop}[thm]{Proposition}

\newtheorem{hyp}[thm]{Hypothesis}
\newtheorem{exam}[thm]{Examples}

\def\pf{{\noindent \bf Proof: }}
\begin{document}

\title[Poisson statistics for  singular randomness]{
Poisson statistics for Anderson model with singular randomness}

\author[Dhriti Ranjan Dolai]{Dhriti Ranjan Dolai}
\address{The Institute of Mathematical Sciences,
    4th Cross Street, CIT Campus,
    Taramani, Chennai 600113, India}
\email{dhriti@imsc.res.in}

\author[M Krishna]{M Krishna}
\address{The Institute of Mathematical Sciences,
    4th Cross Street, CIT Campus,
    Taramani, Chennai 600113, India}
\email{krishna@imsc.res.in}

\thanks{\noindent  
MK was partially supported by IMSc Project 12-R\&D-IMS-5.01-0106 \\
and thanks Peter Hislop, Minami and Anish Mallik for discussions on eigenvalue statistics}

\maketitle

\begin{abstract}
In this work we consider the Anderson model on the d-dimensional lattice with the single site potential having singular distribution, mainly $\alpha$-H\"older continuous ones and show that
the eigenvalue statistics is Poisson in the region of exponential localization.  
\end{abstract}

\tableofcontents

\noindent {\bf  AMS 2000 MSC:} 35J10, 81Q10, 35P20

\section{The Model}

We consider the Anderson model, namely the operators 
\begin{equation}\label{model}
H^\omega = H_0 + \sum_{j \in \ZZ^d} \omega_j P_j, ~~ \omega \in \Omega,
\end{equation}
on $\ell^2(\ZZ^d)$ with $P_j$ the orthogonal projection onto $\ell^2(\{j\})$.
We take $H_0$ to commute with translations on $\ZZ^d$.  Typically we have
$H_0 = \Delta$, which is the discrete Laplacian with diagonal part dropped.
$H_0 = 0$ is also included in the model.
  
We consider a cube $\Lambda_L$ of side length $L$ and cover the cube with smaller disjoint
cubes $C_p$ of side length $l_L$, so that $\Lambda_L = \cup_{p=1}^{N_L^d} C_p$.  Given these we consider
the matrices $H_{\Lambda_L}^\omega$ and $H_{C_p}^\omega$ 
obtained by compressing the operator $H^\omega$ to the finite dimensional
subspaces $\ell^2(\Lambda_L)$ and $\ell^2(C_p)$ respectively.  

An explicit collection of such cubes $C_p$ is obtained by
we dividing $(-L-1,L]^d$ into $N_L^d$ equal cubes $C_p^*$ of the form $(c,d]^d$
for $p=1,2,\cdots,N_L^d$, with side length $\frac{2L+1}{N_L}$ and defining
\begin{equation}\label{cpcubes}
C_p=C_p^*\cap\mathbb{Z}^d, ~~ 
int(C_p)=\{x\in C_p : dist(x,\partial C_p)>l_L\},
\end{equation}
where $\{N_L\}_L$ and $\{l_L\}_L$ are both increasing sequences of integers, which
will be specified latter. For any cube $B\subset\mathbb{Z}^d$ the boundary of $B$ is denoted by 
$\partial B$ and is defined by
$$
\partial B=\{x\in B : \exists~ x'\in B^c~such ~that~|x-x'|=1\}.
$$
\begin{hyp}\label{hyp}
We assume that the single site distribution $\mu$ is uniformly 
$\alpha$-H\"older continuous for $0 < \alpha \leq 1$.  
\end{hyp}
We look in the region of exponential localization to study eigenvalue statistics.  
In this context by 'exponential localization' we mean that some appropriate
fractional moments of the Green functions associated with the operators are exponentially decaying
as shown by Aizenman-Molchanov \cite{AM} in the strong disorder case.  
For simplicity we assume complete localization and state their results as a hypothesis for the present work.  

\begin{hyp}\label{hyp2}
Let $\Lambda \subseteq \ZZ^d$ be any large cube, then the inequality
$$
\sup_{z \in \CC^+, Re(z) \in [a, b] } \EE \left( |\langle \delta_n, (H^\omega_\Lambda - z)^{-1} \delta_m\rangle|^s \right)  \leq C 
e^{- \gamma |n-m|}, 
$$
is valid for some $ s >0 $ and $0 < \gamma < \infty$. 
\end{hyp}
Henceforth the energy $E \in \sigma(H^\omega)$ appearing all the quantities below is assumed
to lie in $[a, b]$ occurring in the above hypothesis.  
 
Using the constant $\gamma$ in the above hypothesis, we specify the 
numbers $N_L, l_L$ of equation (\ref{cpcubes}).  We fix an $0 < \epsilon <1$ and choose
\begin{equation}\label{nl-ll}
N_L = O(L^{1-\epsilon}), ~~ l_L = \frac{5d}{\alpha \gamma} \ln(2L+1).
\end{equation}

Given the operators in equation (\ref{model}) satisfying the above hypotheses
 and an $\alpha$, $0 < \alpha \leq 1$, we define the following random measures on $\RR$,
where we set $ \beta_L = (2L+1)^{\frac{d}{\alpha}}$.
\begin{equation}\label{process}
\begin{split}
\xi_{L,E}^\omega(I) &= Tr(\chi_{\Lambda_L} E_{\beta_L(H_L^\omega - E)}(I))\\
\eta_{p,E}^\omega(I) &= Tr(\chi_{\Lambda_L} E_{\beta_L(H_{C_p}^\omega - E)}(I)), ~ p=1,\dots, N_L^d\\
\zeta_{L,E}^\omega(I) &= Tr(\chi_{\Lambda_L} E_{\beta_L(H^\omega - E)}(I)).
\end{split}
\end{equation}

The eigenvalue statistics was studied for random Schr\"odinger operators by Molchanov \cite{MSA}
followed by Minami \cite{NM}, who obtained an estimate for ensuring that the L\'evy measure
of the eigenvalue point process is degenerate. Eigenvalue statistics was studied by 
Germinet-Klopp \cite{GK}, Aizenman-Warzel \cite{AW} for the canopy graph.  Combes-Germinet-Klein
\cite{JFA} obtained the Minami estimate in a more transparent form while extending the original
estimate to more general single site distributions and obtaining estimates on probabilities
associated with existence of multiple eigenvalues in an interval. The statistics results 
were extended to include localization centers by Nakano-Killip \cite{NK} and for the Schr\"odinger 
case Nakano \cite{Nak} showed infinite divisibility of limiting point processes.  In the case
of Anderson model with higher rank random potentials Tautenhan-Veseli\'c \cite{TV} obtained
the Minami estimate leading to the Poisson statistics.  Recently Hislop-Krishna \cite{HK} 
considered higher rank random potentials for Anderson models and showed the eigenvalue 
statistics to be compound Poisson.  They also showed that in all the cases of random
Schr\"odinger and Anderson models with higher rank i.i.d random potentials, 
the Wegner estimate and complete localization leads to a compound Poisson eigenvalue statistics in general.  

Level repulsion was proved for a class of Anderson models with decaying randomness by 
Dolai-Krishna \cite{DK}, for a class of Schr\"odinger operators in one dimension Kotani-Nakano
\cite{KotNak} obtained $\beta$-ensemble governing the statistics and for localization centers level 
repulsion was shown in the Anderson model by Nakano \cite{Nak2}. 

All these works, except that of Combes-Germinet-Klein \cite{JFA}, assumed that the single site
distributions have an absolutely continuous bounded density, which amounts to taking
$\alpha =1$.   

We show here that by changing the scale appropriately we can include more singular single site 
distributions.  However this comes at a price.  In view of the subtleties 
involved with
singular measures, in particular the absence of a de la Val\'ee Poussin type
theorem, the results become weaker.

It is well known that once we have the Wegner estimate, limit points of 
the point processes in equation (\ref{process}) are also point processes.  
We state this fact as a theorem below.  The proof involves showing 
tightness of the family of measures $\xi_{L,E}^\omega$ and proof is given as part of 
Hislop-Krishna \cite[Proposition 4.1]{HK}, so we state it without proof.
\begin{thm}\label{main0}
Consider $H^\omega$ as in equation (\ref{model}) satisfying Hypothesis \ref{hyp} with 
$0 < \alpha \leq 1$.  Then every limit point, in the 
sense of distributions, of $\xi_{L,E}^\omega$, is a point process.
\end{thm}

The main question is then to determine the nature of the limiting point 
processes which
we do below for $\alpha$-H\"older continuous measures.

To the best of our knowledge ours is 
the first instance where singular single site distributions are allowed to obtain 
eigenvalue statistics.

We use the same symbol $\nn$ for the IDS and the measure associated with it,
and it should be clear from the context the object used.

We define the $\alpha$-derivative and the $\alpha$-upper derivatives
of the integrated density of states $\nn$ of our model, by 
$$
d_{\nn}^\alpha(E) = \lim_{\epsilon \rightarrow 0} \frac{\nn((E-\epsilon, E+\epsilon))}{(2\epsilon)^\alpha}, 
D_{\nn}^\alpha(E) = \limsup_{\epsilon \rightarrow 0} \frac{\nn((E-\epsilon, E+\epsilon))}{(2\epsilon)^\alpha}. ~~ 
$$ 

We define the measures $\ml_\alpha$ as
\begin{equation}\label{elalpha}
\ml_\alpha(I) = \alpha 2^{\alpha -1} \int_I |y|^{\alpha -1} dy
\end{equation}
for any bounded Borel subset $I \subset \RR$.

Our main theorem is then to show that the limiting point processes
give Poisson distributions for a class of intervals.  For technical reasons,
that have to do with the fact that we are dealing with singular measures,
we consider only intervals symmetric about the origin below to obtain the parameters
of the limiting Poisson distributions. 
In view of the fact that for singular measures $\nu$, 
$d_\nu^\alpha(x)$ may not exist for almost all $x$ w.r.t. $\nu$, we
have to deal with upper derivatives in which case we can only talk about
limit points of the random measures we considered above.

\begin{thm}\label{main}
Consider $H^\omega$ as in equation (\ref{model}) satisfying Hypotheses \ref{hyp} with $0 < \alpha \leq 1$
and Hypothesis \ref{hyp2} with $E$ in the region of localization.   
For any bounded open interval $I$, suppose $\gamma_{E,I}$ is non-zero such that 
$$
\gamma_{E,I} = \lim_{n\rightarrow \infty} \EE (\xi_{L_n(I),E}^\omega(I) ). 
$$ 
Then the random variables $\xi_{L_n(I), E}^\omega(I)$ converge in distribution to the Poisson 
random variable with parameter $\gamma_{E, I}$.
\end{thm}

This theorem implies the following.
\begin{cor}\label{cormain}
With the assumptions of Theorem \ref{main}, if  $0 < D_{\nn}^\alpha(E) < \infty$, then
for each bounded open interval $I = -I$, there is a subsequence $L_n(I)$ such that  
$\xi_{L_n(I),E}^\omega(I)$ converges in distribution to a Poisson random variable
with parameter $D_\nn^\alpha(E) \ml_\alpha(I)$. 
\end{cor}

It is interesting to note that the measures $\ml_\alpha$ occur in the theorem of 
Jensen-Krishna, \cite[Theorem 1.3.2]{MM} dealing with continuous wavelet
transforms of measures, where the constants $c_\alpha$
there are integrals of the function $\psi$, the function that generates the "continuous wavelet",
with respect to $\ml_\alpha$.

Since the limsup of a sequence is always a limit point of the sequence, the above theorem shows that when the upper derivative $D_\nn^\alpha(E)$ is positive, there is at least one subsequence of $\xi_{L,E}^\omega(I)$ that
converges in distribution to a Poisson random variable, the parameter
of the corresponding Poisson distribution is then $D_\nn^\alpha(E) \ml_\alpha(I)$.

It was shown by Krishna \cite[Corollary 1.7]{K}, Combes-Germinet-Klein \cite{JFA}, 
( and by Combes-Hislop-Klopp \cite{CHK}, Stollmann \cite{Sto} for continuous models), 
that when the single site distribution is uniformly $\alpha$-H\"older continuous 
(Hypothesis \ref{hyp}), the IDS, 
$\nn$ is uniformly $\alpha$-H\"older continuous.  When $\nn$ is $\alpha$-H\"older continuous
we can use the decomposition of Theorem 69, Rogers \cite{Rog}, with $h(x) = |x|^\alpha$,
to obtain for any bounded Borel subset $E \subset \RR$,
$$
\nn(E) = \int_E f(t) d\mu^h(x) + J(E),  
$$
with $J$ being strongly $h$ continuous.  Accordingly if this $f$ is non-zero a.e. $\mu^h$
then we will have non-zero $D_\nn^\alpha(E)$ for those $E$ for which $f(E)$ is non-zero.

However since the theorems of \cite{K},\cite{CHK},  
may not be optimal, 
in the sense that the $\nn$ may have better modulus of continuity  
in some part of its support, (as seen in Kaminaga-Nakamura-Krishna \cite{KNK}, where
the IDS is analytic in some region of the spectrum even for measures $\mu$ with 
singular component in them) even when $\mu$ is only $\alpha$-H\"older continuous, we 
cannot be sure that the $f$ is indeed non-zero for the given $E$.

We finally note that the subsequences of Theorem \ref{main0} and those in Corollary \ref{cormain}
may be different.

\section{Ideas of Proofs}

The strategy of proof is to first follow the procedure adopted by 
Minami \cite{NM} where one first shows that three
classes of random measures considered above are asymptotically 
essentially the same in the sense that if one of the limits 
below exists then it does for all and they are all the same.
A similar statement holds for any subsequence also.

\begin{prop}\label{prop1}
Consider the processes defined in equation (\ref{process}) associated
with the operators $H^\omega$ satisfying Hypotheses \ref{hyp},\ref{hyp2}. 
Then for any bounded interval $I\subset\mathbb{R}$ we have
\begin{equation}
\label{intensity1}
\lim_{L\to\infty}\zeta^{\omega}_{L,E}(I)= \lim_{L\to\infty}\xi^{\omega}_{L,E}(I)= 
 \lim_{L\to\infty}\sum_p^{N_L^d}\eta^{\omega}_{p,E}(I).
\end{equation}
\end{prop}
In the above the limits are in the sense of convergence in distributions

\begin{cor}\label{cor1}
 For any bounded interval $I\subset\mathbb{R}$ we have
\begin{equation}
\label{intensity2}
\lim_{L\to\infty}\mathbb{E}\big(\zeta^{\omega}_{L,E}(I)\big)= \lim_{L\to\infty}\mathbb{E}\big(\xi^{\omega}_{L,E}(I)\big)= 
 \lim_{L\to\infty}\mathbb{E}\bigg(\sum_p^{N_L^d}\eta^{\omega}_{p,E}(I)\bigg).
\end{equation}
\end{cor}

Once these results are established,
our strategy is to use the the celebrated L\'evy-Khintchine representation theorem
for measures.  The L\'evy-Khintchine theorem (see Theorem 1.2.1, Applebaum \cite{App}) 
says that a measure $\nu$ is infinitely divisible
if and only if its characteristic function $\widehat{\nu}(t)$
 is of the form
$$
e^{iat +b t^2 + \int_{|x|\leq c} (e^{itx} - 1 - itx)dM(x)}
$$
for some $\sigma$-finite measure $M$, which is called the L\'evy measure
associated with $\nu$.  In the case the measure $M$ is finite
we can absorb the linear term into the number $a$ and rewrite this expression
in the form
$$
e^{iat +b t^2 + \int_{|x|\leq c} (e^{itx} - 1)dM(x)}
$$
It turns out that a distribution is Poisson iff
$a=b=0$ and $M$ is supported on $\{1\}$ (notationally such a measure is written
as a positive multiple of $\delta(x-1)dx$ by some authors). 
The parameter of the Poisson distribution is then the number $M(\{1\})$.

We emphasize here that to show Proposition \ref{prop1} it is sufficient
to have exponential localization (in the sense of Aizenman-Molchanov \cite{AM})
and the Wegner estimate for the IDS $\nn$.  The result that 
$M(\RR \setminus \{1\}) = 0$ uses the Minami estimate.

Therefore the idea is to compute the Fourier transforms of the random variables
$$
\sum_{p}\eta_{p,E}^\omega(I)
$$ 
which are a sum of i.i.d random variables and show that the limit of the Fourier
transform has the desired form.  

In view of the Corollary \ref{cor1},  the value of the parameter of the 
Poisson distribution is computed using the fact that the parameter
is the expectation of the Poisson distribution
which in this case is obtained as the limit of 
$\EE\left(\zeta_{L,E}^\omega(I)\right)$ either for the whole sequence
or if the limit does not exits for some subsequences for which it does.

In the context of absolutely continuous single site distributions $\mu$
these limits exist at points in the spectrum where the density of states
exists.  In our context where we are dealing with singular single site
distributions which have no density with respect to the Lebesgue measure
we need to consider derivatives or upper derivatives with respect to
Hausdorff measures to obtain these limit points.

\section{Preliminaries }

It was shown by Krishna \cite[Corollary 1.7]{K} that if the single site distribution $\mu$ is uniformly $\alpha$-H\"{o}lder, 
$0<\alpha\leq 1$ continuous , then the integrated density of states (IDS) is also 
uniformly $\alpha$-H\"{o}lder continuous.\\  We state this fact in the form given
by Combes-Germinet-Klein \cite{JFA}.
Given a probability measure $\mu$ let $S_\mu(s):=\sup_{a\in\mathbb{R}}\mu[a,a+s]$.  Define 
\begin{equation*}
 Q_{\mu}(s):= \left\{
 \begin{array}{rl}
  \|\rho\|_\infty~ s & \text{if } ~\mu ~has~ bounded~ density~ \rho,\\
   8 S_\mu(s) & ~ otherwise.\\
 \end{array} \right.
\end{equation*}
If $\mu$ is a H\"{o}lder continuous with exponent $\alpha\in (0,1]$ then $Q_{\mu}(s)\leq U s^{\alpha}$ for small $s>0$, for some constant $U$.\\
In \cite{JFA}  Combes-Germinet-Klein prove the Wegner estimate and and Minami estimate for more general measure $\mu$ (single site distribution),
we collect their results in the following lemma which immediately gives the following corollary.
The inequality (\ref{Wegner}) is \cite[inequality (2.2)]{JFA}, the inequality 
(\ref{finite volume}) is \cite[Theorem 2.3]{JFA} and the inequality (\ref{Minami}) is \cite[Theorem 2.1]{JFA}, so we omit the proofs.

\begin{lemma}
 \label{Minami-Wegner}
For all bounded interval $I\subset\mathbb{R}$ and  any finite volume $\Lambda\subset\mathbb{Z}^d$, we have
\begin{equation}
 \label{Wegner}
\mathbb{E}\big(\langle \delta_n, E_{H^{\omega}}(I)\delta_n\big) \leq Q_{\mu}(|I|),
\end{equation}
\begin{equation}
 \label{finite volume}
\mathbb{E}\big(Tr(E_{H^\omega_\Lambda}(I))\big)\leq Q_{\mu}(|I|)~|\Lambda|,
\end{equation}
\begin{equation}
 \label{Minami}
\mathbb{E}\bigg(Tr(E_{H^\omega_\Lambda}(I))\big(Tr(E_{H^\omega_\Lambda}(I))-1\big)\bigg)\leq \bigg(Q_{\mu}(|I|)~|\Lambda|\bigg)^2.
\end{equation}
\end{lemma}

\begin{cor}
\label{integrability}
Consider $\nn$, the IDS of the operators $H^\omega$ satisfying Hypothesis \ref{hyp}. Then for  any $\psi \in C_c(\mathbb{R})$ and $n\in\mathbb{Z}^d$ , we have
\begin{equation}
 \int_{\mathbb{R}}\psi(x)d\mathcal{N}(x)=\mathbb{E}\big(\langle \delta_n, \psi(H^{\omega})\delta_n\big)\leq \|\psi\|_{\infty}~Q_\mu(|s_\psi|), ~s_\psi=\mathrm{supp}~\psi.
\end{equation}
\begin{equation}
 \mathbb{E}\big(Tr(f(H^{\omega}_\Lambda))\big)\leq \|\psi\|_{\infty}~Q_\mu(|s_\psi|)~|\Lambda|.
\end{equation}
\end{cor}

Given any measure $\nu$ we denote notationally $\nu(f) = \int f(x) d\nu(x)$ below, where again
the limits are to be understood as in the sense stated for Proposition \ref{prop1}.

\begin{thm}\label{thm2}
Let $H^\omega$ satisfy the Hypotheses \ref{hyp}, \ref{hyp2} and let $E \in \sigma_p(H^\omega)$.
Then for each $\psi \in C_c(\RR)$, we have 
\begin{equation}
\label{vagueconv}
\lim_{L\to\infty}\zeta^{\omega}_{L,E}(\psi)= \lim_{L\to\infty}\xi^{\omega}_{L,E}(\psi)= 
 \lim_{L\to\infty}\sum_p^{N_L^d}\eta^{\omega}_{p,E}(\psi),
\end{equation}
with convergence in the sense of distributions.
\end{thm}
\pf By general theory (see Kallenberg \cite[Theorem 4.5]{kal}), the theorem follows if we show that
\begin{equation}
\label{equ1}
 \lim_{L\to\infty}\bigg|\mathbb{E}\bigg[e^{-\xi_{L,E}^\omega(\psi)}-e^{-\sum_p^{N_L^d}\eta_{p,E}^\omega(\psi)}\bigg]\bigg|=0,
\end{equation}
\begin{equation}
\label{equ2}
 \lim_{L\to\infty}\bigg|\mathbb{E}\bigg[e^{-\zeta_{L,E}^\omega(\psi)}-e^{-\sum_p^{N_L^d}\eta_{p,E}^\omega(\psi)}\bigg]\bigg|=0.
\end{equation}
Since the set of function $\phi_z(x)=Im\frac{1}{x-z},~z\in\mathbb{C}^+$ are dense in $C_c(\mathbb{R})$ it is sufficient to verify (\ref{equ1}) for such function.\\
For $n\in int(C_p)$ and $z\in \mathbb{C}^+$ we have the well known perturbation formula, using the resolvent estimate,
\begin{equation}
 \label{pertur formula}
G^{\Lambda_L}(z;n,n)=G^{C_p}(z;n,n)+\sum_{(m,k)\in\partial C_p}G^{C_p}(z;n,m)G^{\Lambda_L}(z;k,n)
\end{equation}
where $(m,k)\in\partial C_p$ means $m\in C_p$, $k\in \mathbb{Z}^d\setminus C_p$ such that $|m-k|=1$.\\
Denote $z_L=E+\beta_L^{-1}z$ then we have, proceeding as in the proof by Minami \cite{NM},
\begin{align}
\label{estima}
& \bigg|\xi_{L,E}^\omega(\phi_z)-\sum_p^{N_L^d}\eta_{p,E}^\omega(\phi_z)\bigg|\\ 
\nonumber &=
\frac{1}{\beta_{L}}\bigg| TrIm G^{\Lambda_L}(z_L)-\sum_{p} TrImG^{C_{p}}(z_L)\bigg|\\
&\leq\frac{1}{\beta_{L}} \sum_{p} \sum_{n\in C_{p}\backslash int(C_{p})}\big\{ImG^{C_{p}}(z_L;n,n)
 +ImG^{\Lambda_L}(z_L;n,n)\big\}\nonumber\\
 &\qquad  +\frac{1}{\beta_{L}}\sum_{p}\sum_{n\in int(C_{p})}\sum_{(m,k)\in \partial C_{p}}|G^{C_{p}}(z_L;n,m)|
|G^{\Lambda_L}(z_L;k,n)| \nonumber\\
                        &=A_{L}+B_{L}\nonumber.
\end{align}
>From Combes-Germinet-Klein \cite[A.9]{JFA} we have for given $k>0$
\begin{equation}
 \label{diagonal}
Imz~\mathbb{E}(G^\Lambda(z;n,n))\leq \pi\bigg(1+\frac{k}{2}\bigg)S_\mu\bigg(\frac{2~Imz}{k}\bigg).
\end{equation}
Since $Imz_L=\beta_L^{-1}Imz$ with $Im z>0$ so using the $\alpha$-H\"{o}lder continuity of $\mu$ we get 
\begin{align}
 \label{boundary}
\frac{1}{\beta_L}\mathbb{E}\bigg(G^{\Lambda}(z_L;n,n)\bigg) &\leq \frac{1}{Imz}~\pi\bigg(1+\frac{k}{2}\bigg)S_\mu\bigg(\frac{2~Imz_L}{k}\bigg),~~\Lambda=C_p,~\Lambda_L\\
&\leq C ~\bigg(\frac{2~\beta_L^{-1}~Imz}{k}\bigg)^{\alpha}\nonumber\\
&\leq C ~(2L+1)^{-d},~~(\mathrm{since} ~~ \beta_L=(2L+1)^{d/\alpha})\nonumber.
\end{align}
>From the inequality (\ref{estima}) and above we get
\begin{align}
 \label{estimation of A}
\mathbb{E}(A_L) &\leq  C ~(2L+1)^{-d}N_L^d\bigg(\frac{2L+1}{N_L}\bigg)^{d-1}l_L
=O\big(L^{-\epsilon}~ln(L)\big),
\end{align}
in view of our choices for $N_L, l_L$ in equation (\ref{nl-ll}).

On the other hand the term $B_L$ is split as  
\begin{align}
\label{B}
B_L &=\frac{1}{\beta_{L}}\sum_{p}\sum_{n\in int(C_{p})}\sum_{(m,k)\in \partial C_{p}}|G^{C_{p}}(z_L;n,m)|
|G^{\Lambda_L}(z_L;k,n)| \\
   &=\frac{1}{\beta_{L}}\sum_{p}\sum_{n\in int(C_{p})}\sum_{(m,k)\in \partial C_{p}}|G^{\Lambda_L}(z_L;k,n)|^s
|G^{\Lambda_L}(z_L;k,n)|^{1-s}|G^{C_{p}}(z_L;n,m)| \nonumber
\end{align}

Then using the fact that $(m,k)\in \partial C_{p}$ and $n\in int(C_p)$ so that $|n-k|>l_L$ 
for large enough $L$,  the Hypothesis \ref{hyp2} with the number $s$ chosen from
there, to estimate,   
$$
|G^{\Lambda_L}(z_L;k,n)|^{1-s}\leq \frac{1}{|Imz_L|^{1-s}}~\mathrm{and}~
|G^{C_{p}}(z_L;n,m)|\leq \frac{1}{|Imz_L|},
$$
we obtain the following bound from taking expectation in the equality  (\ref{B}).
\begin{equation}
 \label{estimation of B}
\mathbb{E}(B_L)\leq \frac{1}{\beta_L|Imz_L|^{2-s}}N_L^d\bigg(\frac{2L+1}{N_L}\bigg)^d\bigg(\frac{2L+1}{N_L}\bigg)^{d-1}
l_Le^{-\gamma l_L}.
\end{equation}
We simplify the right hand side of the above inequality to get 
\begin{align}
 \label{convergence of B}
\mathbb{E}^{\omega}(B_L)&\leq \frac{1}{\beta_L|Imz_L|^{2-s}}N_L^d\bigg(\frac{2L+1}{N_L}\bigg)^d\bigg(\frac{2L+1}{N_L}\bigg)^{d-1}l_Le^{-rl_L} \\
&=O\big(L^{(1-s)\frac{d}{\alpha}} N_L^d \big(\frac{L}{N_L}\big)^{2d-1} e^{-\gamma \frac{5d}{\gamma \alpha} \ln(2L+1)} \ln(L) \big) \\
&=O\big(L^{(\frac{d}{\alpha}(1-s) + (2d-1) + (d-1)(1-\epsilon) - \frac{5d}{\alpha})} \ln(L)  \big) \\
&= O(L^{-1}), 
\end{align}
since $\alpha \leq 1$.  
In particular we have from (\ref{estimation of A}) and (\ref{convergence of B})
\begin{equation}
 \label{mean converge }
\mathbb{E}(A_L+B_L)\rightarrow 0~~as~~L\to\infty.
\end{equation}
Finally the inequality
$|e^{-x}-e^{-y}|\leq |x-y|~for~x,y>0$ together with the bound 
(\ref{estima}) and above convergence gives the required vanishing of the limits
equation (\ref{equ1}).\\
Again using the resolvent equation for $G(z;n,n)=\langle\delta_n,(H^{\omega}-z)^{-1}\delta_n\rangle$
and the equality
$$
\langle \delta_n, (\oplus H^\omega_{C_p} \oplus H^\omega_{\Lambda_L^c} - z)^{-1} \delta_n\rangle
= \langle \delta_n, (\oplus H^\omega_{C_p}  - z)^{-1} \delta_n\rangle
$$
valid for each $n \in \Lambda_L$, gives us  the relation
\begin{equation*}
G(z;n,n)=G^{C_p}(z;n,n)+\sum_{(m,k)\in\partial C_p}G^{C_p}(z;n,m)G(z;k,n), ~~ n \in C_p, 
\end{equation*}
for each $p=1,\dots, N_L^d$.  The convergence in equation (\ref{equ2}) is then obtained  by
essentially repeating the argument above. \qed

\section{Proof of the main Theorem }

We first prove the Theorem \ref{main}.

In the argument below we consider a subsequence $L_n$ which converges to
the limsup in equation (\ref{gamma}) and use  
Proposition \ref{prop1} to only consider $\xi_{L,E}$ instead of $\zeta_{L,E}$.
We will show that 
$$
\lim_{n\rightarrow \infty} \EE\big(e^{it\xi_{L_n,E}^\omega(I)}\big) = e^{(e^{it} -1)\gamma_{E,I}}. 
$$
This will then show, by L\'evy-Khintchine theorem that $\xi_{L_n,E}^\omega(I)$ converge
in distribution to the Poisson random variable with parameter  $\gamma_{E,I}$.
Since the convergence in distribution for a sequence of random variables is equivalent 
to the convergence of their Fourier transforms point wise combined with Theorem \ref{thm2}, 
it is enough to look at the limit with
$\xi_{L_n,E}^\omega(I)$ replaced by $\sum_{p=1}^{N_{L_n}^d} \eta_{L_n,E}^\omega(I)$.

We first note that from equation (\ref{intensity2}) we have,
\begin{align}
\label{derivative}
 \lim_{n\to\infty}\mathbb{E}\bigg(\sum_p^{N_{L_n}^d}\eta^{\omega}_{p,E}(I)\bigg)
 =\lim_{n\to\infty}\mathbb{E}\big(\zeta^{\omega}_{L_n,E}(I)\big) 
=\gamma_{E,I}\nonumber\\
\end{align}
We now compute the limits of Fourier transforms 
\begin{align}
\label{prod}
\lim_{n\rightarrow \infty} \EE\big(e^{it\xi_{L_n,E}^\omega(I)}\big) &= \lim_{n\rightarrow \infty} \mathbb{E}\bigg(e^{it\sum_{p=1}^{N_{L_n}^d}\eta^\omega_{p,E}(I)}\bigg)\\
&=\lim_{n\rightarrow \infty} \prod_{p=1}^{N_{L_n}^d}\mathbb{E}\big(e^{it\eta^\omega_{p,E}(I))}\big)\nonumber\\
&=\lim_{n\rightarrow \infty} \bigg[\mathbb{E}\big(e^{it\eta^\omega_{p,E}(I))}\big)\bigg]^{N_{L_n}^d}\nonumber
\end{align}
Now for $p=1,\dots, N_{L_n}^d$,
\begin{align}
\label{cal}
 \mathbb{E}\big(e^{it\eta^\omega_{p,E}(I))}\big)&=\sum_{m=0}^\infty e^{itm}\mathbb{P}\big(\eta^\omega_{p,E}(I))=m\big)\\
&=1+\mathbb{E}(\eta^\omega_{p,E}(I))[e^{it}-1]+R_{L_n}\nonumber
\end{align}
where $R_{L_n}$ is given by
\begin{align*}
 R_{L_n} &=\sum_{m=0}^\infty e^{itm}\mathbb{P}\big(\eta_{p,E}^{\omega}(I))=m\big)-1-\mathbb{E}(\eta_{p,E}^{\omega}(I))[e^{it}-1]\\
&=\sum_{m=0}^\infty e^{itm}\mathbb{P}\big(\eta_{p,E}^{\omega}(I))=m\big)-\sum_{m=0}^\infty \mathbb{P}\big(\eta_{p,E}^{\omega}(I)=m) \\ 
&  ~~ -[e^{it}-1]
\sum_{m=0}^\infty m\mathbb{P}\big(\eta_{p,E}^{\omega}(I))=m\big)\\
&=\sum_{m=2}^\infty \big(e^{itm}-me^{it}+m-1\big)\mathbb{P}\big(\eta_{p,E}^{\omega}(I))=m\big).
\end{align*}
Then using the inequality 
$$
|e^{itm} - me^{it} + m -1| \leq (m+1) + (m-1) \leq 2m, ~~ \mathrm{when} ~~ m \geq 2
$$
and setting  $J_{L,E}=E+\beta_{L}^{-1}I$ we get,
\begin{align*}
|R_{L_n}| &\leq\sum_{m=2}^\infty\big(|e^{itm}-me^{it}|+(m-1)\big)~~\mathbb{P}\big(\eta_{p,E}^{\omega}(I))=m\big)\\
 &\leq\sum_{m=2}^\infty\big((m+1)+(m-1)\big)~~\mathbb{P}\big(\eta_{p,E}^{\omega}(I))=m\big)\\
&\leq 2\sum_{m=2}^\infty m\mathbb{P}\big(\eta_{p,E}^{\omega}(I))=m\big)\\
&\leq 2\sum_{m=2}^\infty m(m-1)\mathbb{P}\big(\eta_{p,E}^{\omega}(I))=m\big)\\
&\leq 2\mathbb{E}\bigg(Tr(E_{H^\omega_{C_p}}(J_{L_n,E})\big(Tr(E_{H^\omega_{C_p}}(J_{L_n,E})-1\big)\bigg),
\end{align*}
Now  from the Minami estimate of Lemma \ref{Minami-Wegner}
inequality (\ref{Minami}) we have
\begin{align}
 \label{double point}
& N_{L_n}^d~\mathbb{E}\bigg(Tr(E_{H^\omega_{C_p}}(J_{L_n,E})\big(Tr(E_{H^\omega_{C_p}}(J_{L_n,E})-1\big)\bigg) \\ &\leq N_{L_n}^d\bigg(Q_{\mu}(|J_{L_n,E}|)~|C_p|\bigg)^2 \nonumber \\
&\leq N_{L_n}^d\big(|J_{L_n,E}|^\alpha |C_p|\big)^2\nonumber\\
&=O\bigg(\beta_{L_n}^{-2\alpha}N_{L_n}^d\bigg(\frac{2L_n+1}{N_{L_n}}\bigg)^{2d}\bigg).\nonumber
\end{align}
The above calculation together with (\ref{double point}) estimate will give
$$
N_{L_n}^dR_{L_n}\rightarrow 0~~as~~n\to\infty.
$$
>From the above computation we get
$$
N_{L_n}^d\bigg[\mathbb{E}(\eta_{p,E}^{\omega}(I))[e^{it}-1]+R_{L_n}\bigg]\xrightarrow{n\to\infty}\gamma_{E,I} [e^{it}-1].
$$
We use the equations (\ref{prod}) and (\ref{cal}) to obtain the equality 
\begin{align}
 \mathbb{E}\big(e^{it\xi_{L_n,E}^\omega(I)}\big)&=\bigg[1+\frac{N_{L_n}^d\big[\mathbb{E}(\eta_{p,E}^{\omega}(I))[e^{it}-1]+R_{L_n}\big]}{N_{L_n}^d}\bigg]^{N_{L_n}^d}.
\end{align}
which combined with the convergence of $\big(1+\frac{z_n}{n}\big)^n$ to $e^z$, 
whenever $z_n\to z$ as $n\to\infty$ gives us finally 
the limit 
\begin{equation*}
 \mathbb{E}\big(e^{it\xi_{L_n,E}^\omega(I)}\big)\xrightarrow{n\to\infty}e^{\gamma_{E,I}(e^{it}-1)}.
\end{equation*}
\qed

\leftline{\bf Proof of Corollary \ref{cormain}:}

We first note that, if we denote 
$$
D_\mu^\alpha(E) < \infty ~ \mathrm{iff} ~ \limsup \frac{\nn(E +\epsilon I)}{(2\epsilon)^\alpha} < \infty, ~~ \mathrm{for ~ all ~ bounded ~~ symmetric ~~ intervals } ~ I.
$$

We will show that 
$$
\limsup \frac{\big(\beta_L^{\alpha} \EE \big( \langle \delta_0, E_{H^\omega}(E + \beta_L^{-1} I) \delta_0\rangle \big)}{|I|^\alpha} \geq \frac{1}{2^d} D_\nn^\alpha(E).
$$
Then by the assumption of theorem the right hand side is positive, so a limit point
of $\EE(\zeta_L(I))$ is positive.  
We recall that $\beta_L^\alpha = (2L+1)^d$.  Let $I$ be a bounded open interval and choose 
$\epsilon \in (\beta_{L+1}^{-1}, \beta_L^{-1}]$,then we have
$$
E+\beta_L^{-1}I\supseteq E+\epsilon I~~\mathrm{and} ~~ \mathcal{N}\big(E+\beta_L^{-1}I\big)\geq \mathcal{N}\big( E+\epsilon I\big).
$$
Therefore we have, since $\beta_{L+1}^\alpha \epsilon^\alpha\ge 1$,
\begin{align}
\frac{\beta_L^\alpha\mathcal{N}\big(E+\beta_L^{-1}I\big)}{|I|^\alpha} &\ge \bigg(\frac{\beta_L}{\beta_{L+1}}\bigg)^\alpha\frac{\mathcal{N}\big( E+\epsilon I\big)}{(\epsilon |I|)^\alpha}~~~
\\
&\ge  \bigg(\frac{\beta_L}{\beta_{L+1}}\bigg)^\alpha\frac{\mathcal{N}(E-c\epsilon,E+c\epsilon)}{(\epsilon |I|)^\alpha}\nonumber,
\end{align}
These inequalities imply that
\begin{align}
\sup_{L\ge M}\frac{\beta_L^\alpha \mathcal{N}\big(E+\beta_L^{-1}I\big)}{|I|^\alpha} &\ge \bigg(\frac{1}{1+\frac{2}{2M+1}}\bigg)^d \sup_{\epsilon \in (\beta_{L+1}^{-1}, \beta_L^{-1}],~L\ge M}
\frac{\mathcal{N}(E + \epsilon I)}{(\epsilon |I|)^\alpha} \\
& \ge  
 \bigg(\frac{1}{1+\frac{2}{2M+1}}\bigg)^d \sup_{\epsilon \in (0, \beta_M^{-1}]} \frac{\mathcal{N}(E + \epsilon I)}{(\epsilon |I|)^\alpha},
\end{align}
where  we used the facts that
$$
\bigcup_{L\geq M} \big(\beta_{L+1}^{-1}, \beta_{L}^{-1}\big] = \big(0, \beta_M^{-1}\big] ~~ \mathrm{and} ~~ \frac{\beta_L}{\beta_{L+1}} \geq \frac{1}{1+\frac{2}{2M+1}}, ~~\mathrm{for} ~~ L \geq M. 
$$
We now let $M\to\infty$ in both side of above then from the definition of limsup we get
\begin{equation}\label{lower}
 \varlimsup_{L\to\infty}\frac{\beta_L^\alpha \mathcal{N}\big(E+\beta_L^{-1}I\big)}{|I|^\alpha}\ge  D^\alpha_{\mathcal{N}}(E).
\end{equation}
Similarly starting with  $\epsilon \in (\beta_{L+1}^{-1}, \beta_L^{-1}]$
we get the inequality
$$
\frac{\beta_{L+1}^\alpha\mathcal{N}\big(E+\beta_{L+1}^{-1}I\big)}{|I|^\alpha} 
\le \bigg(\frac{\beta_{L+1}}{\beta_{L}}\bigg)^\alpha\frac{\mathcal{N}\big( E+\epsilon I\big)}{(\epsilon |I|)^\alpha}~~~
$$
and proceed as in the above argument, with upper bounds now, to get
\begin{equation}\label{upper}
 \varlimsup_{L\to\infty}\frac{\beta_L^\alpha \mathcal{N}\big(E+\beta_L^{-1}I\big)}{|I|^\alpha}\le  D^\alpha_{\mathcal{N}}(E).
\end{equation}
Putting the inequalities (\ref{lower}) and (\ref{upper}) we get
$$
\varlimsup_{L\to\infty}\frac{\beta_L^\alpha \mathcal{N}\big(E+\beta_L^{-1}I\big)}{|I|^\alpha} =  D^\alpha_{\mathcal{N}}(E).
$$
The above inequality shows that, noting again that $\beta_L^\alpha = (2L+1)^d$,
\begin{align}\label{gamma}
 \gamma_{E,I} 
&= \varlimsup_{L\to \infty} \EE\big(\zeta_{L,E}^\omega(I)\big)  \\ 
&= \varlimsup_{L\to \infty} \EE\big(\sum_{n\in \Lambda_L} \langle\delta_n, E_{H^\omega}(I) \delta_n\rangle \big) \nonumber \\ 
&=  \varlimsup_{L\to\infty}\beta_L^\alpha \mathcal{N}\big(E+\beta_L^{-1}I\big) \nonumber  \\
&= D^\alpha_{\mathcal{N}}(E) |I|^\alpha  =  D^\alpha_{\mathcal{N}}(E) \ml_\alpha(I) \nonumber  ,
\end{align}
where to pass to the third line we used the fact that 
$\EE\big(\langle\delta_n, E_{H^\omega}(I) \delta_n\rangle\big)$ does not depend on $n$.  Since the limsup above is a limit point of the sequence considered, we have the corollary.
\qed

\section{Example}

\begin{exam}
We now give an example of random operators that have singular density of states and for which 
the local eigenvalue statistics is Poisson.  We note  while this example may appear
trivial, it is one for which none of the existing theorems can show Poisson eigenvalue statistics.

Consider the operator
$$
H^\omega = \sum_{n\in \ZZ^d} \omega_n P_n
$$
$P_n$ is projection onto $\ell^2(\{n\})$ as in the model (\ref{model}) with 
$\{\omega_n\}$ i.i.d random variable distributed by a measure $\mu$.  Then the IDS agrees
with the distribution of the measure $\mu$, so if we choose a singular $\alpha$-continuous
measure $\mu$ (such as the Cantor measure, for which $\alpha = \log(2)/\log(3)$), then
the conditions of our theorem are valid for $H_0= 0$ (which is in some sense infinite disorder
limit of the large disorder Anderson model).  

Therefore Poisson eigenvalue statistics holds for points in the spectrum.
\end{exam}

\end{document}